\documentclass[11pt,english,a4paper]{article}

\usepackage{babel}
\usepackage{amsmath,amsfonts,amssymb,amsthm}
\usepackage{color}

\theoremstyle{plain}% default
\newtheorem{thm}{Theorem}% [section]

\newtheorem{prop}{Proposition}

\theoremstyle{definition}

\newtheorem{case}{Case}

\theoremstyle{remark}

\newcommand{\Z}{\mathbb{Z}}
\newcommand{\N}{\mathbb{N}}

\newcommand{\T}{\mathbb{T}}

\newcommand{\F}{\mathbb{F}}

\newcommand{\hi}{\mathcal{H}}

\newcommand{\m}{\text{-}}
\newcommand{\inv}[1]{#1^{\m 1}}

\newcommand{\ospec}[1]{{\widehat{#1}}^{\, o}}

\DeclareMathOperator*{\aut}{Aut}

\def\H{\mathcal{H}}
\def\K{\mathcal{K}}

\begin{document}

\title{Primitivity of some full group C$^*$-algebras}

\author{Erik B\'edos$^{*}$,
Tron {\AA}. Omland$^{*}$ \\
\\
\\}

\date{{\it March 29, 2010}}
\maketitle
\renewcommand{\sectionmark}[1]{}

\begin{abstract}
We show that the full group C$^*$-algebra of the free product of two nontrivial countable amenable discrete groups, where at least one of them has more than two elements, is primitive. We also show that in many cases, this C$^*$-algebra is antiliminary and has an   
uncountable family of pairwise inequivalent, faithful irreducible representations.


\end{abstract}

\vfill
\thanks{\noindent $^*$ partially supported by the Norwegian Research
Council.\\

\par \par
 
 \noindent E.B.{'}s address: Institute of Mathematics, University of
Oslo, P.B. 1053 Blindern, 0316 Oslo, Norway. E-mail : bedos@math.uio.no. \\

\noindent 
T.O.{'}s address: Department of Mathematical Sciences, NTNU,
7491 Trondheim, Norway. E-mail : tronanen@math.ntnu.no.\par}

\begin{flushleft}

\section{Introduction} 

Let $G$ denote a countable discrete group. It is known that $C^{*}(G)$, the full group C$^*$-algebra of $G$, is primitive in a number of cases \cite{Y, Cho, Pa, Kha, Mur, BO}. Especially, this is true for many groups which have a free product decomposition satisfying various  conditions: see \cite{Kha, Mur, BO}.  These results suggest that $C^{*}(G)$ should be primitive whenever $G$ is the free product  of two nontrivial countable discrete groups $G_1$ and  $G_2$, where  at least one of them has more than two elements. In this note, we show that this is indeed the case when both $G_1$ and $G_2$ are also assumed to be amenable. 

This applies for example when $G_1$ and $G_2$ are both finite with $|G_1| \geq 2$ and $|G_2| \geq 3$. This case is not covered by any of the papers cited above, except when $G_1 = \Z_2$ and $G_2= \Z_3$, i.e.\  $G$ is the modular group $PSL(2, \Z)$, for which primitivity of $C^*(G)$ was shown in \cite{BO}.  The reader should consult \cite{Mur} and \cite{BO} for more information around the problem of determining when the full group C$^*$-algebra of a countable discrete group is primitive.

\bigskip Our proof will rely on the following result from \cite{BO}:
\begin{thm}\label{theorem 1}
Assume that a group $G$ has a normal subgroup $H$ such that
\begin{itemize}
\item[i)] $C\sp*(H)$ is primitive,
\item[ii)] $K=G/H$ is amenable,
\item[iii)] the natural action of $K$ on $\ospec{C\sp*(H)}$  has a free point.
\end{itemize}
Then $C\sp*(G)$ is primitive.
\end{thm}
We recall here what condition iii) means. 
Set  $A=C^*(H)$. Then $ \ospec{A}=\{[\pi]\in\widehat{A}\mid\pi\text{ is faithful}\}$ is nonempty since $A$ is assumed to be primitive. The natural action of $K=G/H$ on $\ospec{A}$ is defined as follows.

\medskip Let $n\colon K\to G$ be a normalized section for the canonical homomorphism $p$ from $G$ onto $K$. Let $\alpha\colon K\to\aut{(A)}$ and $u\colon K\times K\to A$ be given by
\begin{equation*}
\begin{split}
\alpha_{k}\big(i_{H}(h)\big)= i_{H}\big(n(k)\, h \, \inv{n(k)}\big),\quad & k\in K, h\in H,\\
u(k,l)=i_{H}\big(n(k)\, n(l)\, \inv{n(kl)}\big),\quad & k,l\in K,
\end{split}
\end{equation*}
where $i_{H}$ denotes the canonical injection of $H$ into $A$.

\medskip Then $(\alpha,u)$ is a twisted action of $K$ on $A$ (cf.\ \cite{PR}), which 
induces an action of $K$ on $\ospec{A}$ given by
\begin{equation*}
k\cdot [\pi]= [\pi\circ\alpha_{\inv{k}}]\,  .
\end{equation*}
This action is independent of the choice of normalized section for $p$ and called the natural action of $K$ on $\ospec{A}$. Finally, we recall that $[\pi] \in \ospec{A}$ is  a free point for this action whenever we have  $k\cdot [\pi]\neq [\pi] $ for all $k\in K,k\neq e$.

\bigskip Throughout this paper, we let $G_1$ and $G_2$  be two nontrivial  countable discrete groups and assume that at least one of them has more than two elements. Further we let $G=G_1 * G_2$ denote the free product of $G_1$ and $G_2$. It is well known that $G$ is icc and nonamenable. Section 2 is devoted to the proof of our main result in this paper:

\begin{thm}\label{theorem 2}
Assume moreover that $G_1$ and $G_2$ are both amenable. 
Then $C\sp*(G)$ is primitive.
\end{thm}

In the final section, we discuss the problem of deciding when $C^*(G)$ is antiliminary and has an uncountable family of pairwise inequivalent, faithful irreducible representations. 

\medskip As will be evident from its proof, the annoying amenability assumption in Theorem 2
is due to the amenability assumption on $K$  in Theorem 1. Now, if one replaces this assumption on $K$ by requiring that  the twisted action of $K$ on $C^*(H)$ is amenable in the sense that the full and the reduced crossed products of $C^*(H)$ by this action agree, then Theorem 1 still holds. An interesting problem  is whether one can find condition(s) other than the amenability of $K$ ensuring that this more general requirement is satisfied.

\section{Proof of Theorem 2}

We let $e_1$ (resp. $e_2$) denote the unit of $G_1$ (resp. $G_2$)   and set $G'_1=  G_1 \setminus \{e_1\}, \, G'_2=  G_2 \setminus \{e_2\}$.
We let $X\subset G$ denote the set of commutators given by 
\begin{equation*}
X=\{\, [a, b] = a\, b\, \inv{a}\inv{b} \in G \, \mid \,  a\in G'_1 , \, b \in G'_2 \,  \} \, .
\end{equation*}
As is well known (see e.g.\ \cite{Se}),  $X$ is free and generates the kernel $H$ of the canonical homomorphism $p$ from the free product 
$G=G_1*G_2$ onto the direct product $K = G_1\times G_2$. The map $(a,b) \to [a,b]$ is then a bijection between $G'_{1} \times G'_{2}$ and $X$,
 and $H$ is isomorphic to the free group $\F_{\lvert X\rvert}$ with $|X|$ generators. 

\medskip  As 
$|X | = |G'_1| \cdot |G'_2 | \geq  2$ , $A=C^*(H)$ is primitive (cf.\ \cite{Y, Cho}). 

\medskip Further, as  $G_1$ and $G_2$ are both assumed to be amenable, $K$ is amenable.

\bigskip Now let $\pi$ be a faithful irreducible representation of $A$ acting on a (necessarily separable) Hilbert space $\hi_{\pi}$.
For each function $\lambda\, \colon \, X\to\T$, we let $\gamma_{\lambda}$ denote the $*$-automorphism of $A$  determined by
\begin{equation*}
\gamma_{\lambda}(i_H(x))=\lambda(x)i_H(x)\, , \quad x \in X ,
\end{equation*}
and set $\pi_{\lambda}=\pi\circ\gamma_{\lambda}$.
Clearly, each $\pi_{\lambda}$ is also faithful and irreducible, i.e. $[\pi_\lambda] \in \ospec{A}$. 

\bigskip The burden of the proof is to establish the following:

\begin{prop} There exist   $[\pi] \in \ospec{A}$ 
and  $\lambda \, \colon \, X\to\T$ such that $[\pi_{\lambda}]$ is a free point for the natural action of $K$ on $A$.

\end{prop}
Once we have  proven this proposition,  the primitivity of $C^*(G)$  then clearly follows from Theorem 1 and the proof of Theorem 2 will therefore be finished.

\begin{proof}[Proof of  Proposition 1]

\medskip 

As a normalized section $n : K \to G$ for $p$, we choose $$n(a,b) = a\, b\, , \quad a \in G_1, \, b \in G_2\, .$$  

We have to show that some faithful irreducible representation $\pi$ of $A$ and some $\lambda\,  \colon \, X\to\T$ may be chosen so that 
\begin{equation*}
\pi_{\lambda}\circ\alpha_k\not\simeq \pi_{\lambda}
\end{equation*}
for all nontrivial $k\in K$.

\medskip 
Clearly, to show that this condition holds, it suffices to show that for each nontrivial $k\in K$, there exists  some $x \in X$ (depending on $k$) such that
\begin{equation}\label{eq:free}
(\pi_{\lambda}\circ\alpha_k)(i_H(x))\not\simeq \pi_{\lambda}\big(i_H(x)\big)\, .
\end{equation}

To show this, we will use following fact:

\medskip  Assume $x_0 \in X$ is fixed. Then, as follows from Choi's  proof \cite{Cho} (see  \cite[Proof of Theorem 3.2]{Mur}), we may choose a faithful irreducible representation $\pi$ $=\pi_{x_0}$ of $A$ such that for each $x \neq x_0$ in $X$ the unitary operator $\pi(i_H(x))$ is diagonal relative to some orthonormal basis of $\hi_\pi$ (which depends on $x$). We will call such a representation for a {\it Choi representation of} $A$ {\it associated to} $x_0$.  

\medskip Our choice of $x_0$, and thereby of $\pi$ $=\pi_{x_0}$, will depend on the possible existence of elements of order $2$ in $G_1$ or $G_2$.

\medskip We will also use repeatedly the following elementary fact (already used in \cite{Mur} and in \cite{BO}): 

\medskip Assume $\H$ is a separable Hilbert space. Let $U$ and $V$ be unitary operators on $\H$ and assume that $U$ is diagonal relative to some orthonormal basis of $\H$. Then the sets   

\begin{equation*}
\{\mu \in \T \mid \mu \, U  \simeq V \}\ \text{and} \  \{\mu \in \T \mid \mu \,  U \simeq (\mu \, U)^* \}
\end{equation*}
are both countable.

\bigskip Consider some faithful irreducible representation $\pi$ of $A$ and $\lambda\,  \colon \, X\to\T$.

\medskip When  $a\in G'_1 , \, b \in G'_2$, so $[a,b] \in X$, we let  $U(a,b)$ (= $U_\pi(a,b)$) denote the unitary operator on $\hi_{\pi}$ given by $U(a,b)= \pi\big(i_H([a,b])\big)$. Further, we set $\lambda(a,b)=\lambda([a,b])$. Thus we have 
\begin{equation}\label{eq:piab}
\pi_{\lambda}\big(i_H([a,b])\big)=\lambda(a,b)U(a,b)\, .
\end{equation}
Some straightforward calculations give the following identities which we will use in the sequel: 
\begin{equation}\label{eq:ab}
\begin{split}
\pi_{\lambda}\big(\alpha_{(a,b)}(i_H([a^{-1},b^{-1}]))\big)&=\lambda(a,b)U(a,b)\\
\pi_{\lambda}\big(\alpha_{(a,e_2)}(i_H([a^{-1},b]))\big)&=\big(\lambda(a,b)U(a,b)\big)\sp*\\
\pi_{\lambda}\big(\alpha_{(e_1,b)}(i_H([a,b^{-1}]))\big)&=\big(\lambda(a,b)U(a,b)\big)\sp*
\end{split}
\end{equation}

\begin{equation}\label{eq:abc}
\begin{split}
\pi_{\lambda}\big(\alpha_{(a,b)}(i_H([a^{-1}c,b^{-1}]))\big)&=\lambda(a,b)U(a,b)\, \big(\lambda(c,b)U(c,b)\big)\sp*\\
\pi_{\lambda}\big(\alpha_{(a,b)}(i_H([c,b^{-1}]))\big)&=\lambda(a,b)U(a,b)\, \big(\lambda(ac,b)U(ac,b)\big)\sp*
\end{split}
\end{equation}

\begin{equation}\label{eq:aecb}
\begin{split}
\pi_{\lambda}\big(\alpha_{(a,e_2)}(i_H([a^{-1}c, b]))\big)&=\lambda(c,b)U(c,b)\, \big(\lambda(a,b)U(a,b)\big)\sp*\\
\pi_{\lambda}\big(\alpha_{(a,e_2)}(i_H([c,b]))\big)&=\lambda(ac,b)U(ac,b)\, \big(\lambda(a,b)U(a,b)\big)\sp*
\end{split}
\end{equation}
whenever $a \in G'_1$, $ \, b \in G'_2$ and $c \in G'_1\setminus \{a, a^{-1}\}$.

\bigskip

We now show how to pick $\pi$ and $\lambda$ such that \eqref{eq:free} holds.
It turns out that the possible existence of elements of order $2$ in $G_1$ or $G_2$  complicates the argument. 
Set  $\ P = \{ s \in G'_1 \, | \, s^2 \neq e_1\}\ , \ Q = \{ t \in G'_2 \, | \, t^2 \neq e_2\}\, ,$

\medskip and $\ S = G'_1 \setminus P\ , \ T = G'_2 \setminus Q\, .$ Hence 
$$G_1=\{e_1\}\sqcup P \sqcup S \, \, \text{and} \, \, G_2=\{e_2\} \sqcup Q \sqcup T\, .$$

\smallskip We divide our discussion into three separate cases. 

\begin{case}

{\it Both $P$ and $Q$ are nonempty}.

\medskip We  pick  $p_0 \in P\, ,\, q_0\in Q$ and set  $x_0 = [p_0^{-1}, q_0^{-1}]  \in X$. 

\medskip Then we let   $\pi=\pi_{x_0}$ be a Choi representation of $A$ associated to $x_0$, and set $U(a,b) = U_{\pi}(a,b)$ for each $x=[a,b] \in X$. 

\medskip It remains to define $\lambda : X \to \T$ so that \eqref{eq:free} holds for each nontrivial $k \in K$.

\medskip We introduce the following notation. 

\medskip 
Assume that  $a \in G'_1\, , \, b \in G'_2\, ,\, p \in P\, , \, q \in Q\, , s\in S, t \in T$. Then we set  
\begin{equation*}
\begin{split}
\Omega(a,b)&=\{\mu  \in \T \mid \mu\, U(a,b)\simeq U(a^{-1},b^{-1})\}\, ,\\
\Omega_1(p)&=\{\mu \in \T \mid \mu \,U(p,q_0)\simeq U(p^{-1},q_0)\sp*\, \} \, ,\\
\Omega_2(q)&=\{\mu \in \T \mid \mu \,U(p_0,q)\simeq U(p_0, q^{-1})\sp*\, \} , \\
\Omega_1(s)&=\{\mu \in \T \mid \mu \,U(s,q_0)\simeq \big (\mu \, U(s,q_0)\big)\sp*\, \} \, ,\\
\Omega_2(t)&=\{\mu \in \T \mid \mu \,U(p_0,t)\simeq \big(\mu \, U(p_0, t)\big)\sp*\, \} \, . \\
\end{split}
\end{equation*}

\smallskip Note that if $(a, b)\neq (p_0^{-1}, q_0^{-1})$, then $\Omega(a,b)$ is countable (as $U(a,b)$ is then diagonalisable). 

Similarly, $\Omega_1(p), \, \Omega_2(q), \,  \Omega_1(s)$ and $\Omega_2(t)$ are countable.

\medskip To ease our notation, we will define $\lambda$ on $G'_1 \times G'_2$ and identify it with the function on $X$ given by $\lambda ([a,b]) = \lambda(a,b),\,  a \in G'_1,\, b \in G'_2$.  

\medskip We will first define $\lambda$ on  $P \times Q$.

\bigskip 
Let  $P= \sqcup_{i \in I} \, \{p_i, p_i^{\, -1}\}\, , \, Q= \sqcup_{j \in J} \, \{q_j, q_j^{\, -1}\}$ be enumerations of $P$ and $Q$, where the index set $I$ (resp.\ $J$) is 
a (finite or infinite) set of successive integers starting from 0.  

\bigskip For each $ i \in I$ and $j \in J$, 
we set 
$$\lambda(p_i^{-1},q_j)=\lambda(p_i^{-1},q_j^{-1}) = 1 \, .$$

Now let $ i \in I , j \in J$. Using \eqref{eq:piab} and \eqref{eq:ab}, we see that \eqref{eq:free} will hold for
\begin{equation*}
\begin{split}
k=(p_i,q_j^{-1}) \text{ and }k=(p_i^{-1},q_j)\text{ if }\ &\lambda(p_i, q_j^{-1})U(p_i, q_j^{-1})\not\simeq U(p_i^{-1},q_j)\, ; \\
k=(p_i,q_j) \text{ and }\ k= (p_i^{-1},q_j^{-1})\text{ if }\ &\lambda(p_i,q_j)U(p_i,q_j)\not\simeq U(p_i^{-1}, q_j^{-1})\, ; \\
k= (p_i,e_2)\text{ and }k=(p_i^{-1},e_2)\text{ if }\ &\lambda(p_i,q_0)U(p_i,q_0)\not\simeq U(p_i^{-1},q_0)\sp* \, ; \\
k=(e_1,q_j) \text{ and }k=(e_1,q_j^{-1})\text{ if}\ &\lambda(p_0,q_j)U(p_0,q_j)\not\simeq \Big(\lambda(p_0,q_j^{-1})U(p_0,q_j^{-1})\Big)\sp* \, .
\end{split}
\end{equation*}

\bigskip

For each $i \in I$ and $ j\in J$, we therefore pick
\begin{equation*}
\lambda(p_i,q_j^{-1})\in\T\setminus\Omega(p_i,q_j^{-1}).
\end{equation*}

Next, for each $i \in I,\, i\neq 0,$ and $j \in J, \, j\neq 0,$ we pick
\begin{equation*}
\begin{split}
\lambda(p_i,q_j)&\in\T\setminus \Omega(p_i,q_j) \, , \\
\lambda(p_i,q_0)&\in\T\setminus\Big(\Omega(p_i,q_0)\cup\Omega_1(p_i)\Big)\, , \\
\lambda(p_0,q_j)&\in\T\setminus\Big(\Omega(p_0,q_j)\cup\lambda(p_0,q_j^{-1})\Omega_2(q_j)\Big)\, .\\
\end{split}
\end{equation*}
Finally, we pick
$$\lambda(p_0,q_0)\in\T\setminus\Big(\Omega(p_0,q_0)\cup\Omega_1(p_0)\cup\lambda(p_0,q_0^{-1})\Omega_2(q_0)\Big)\, .$$
All these choices are possible as all the involved $\Omega$'s are countable. After having done this,  $\lambda $ is defined on $P\times Q$ and we know that
\eqref{eq:free} will hold for all $k \in (P\times Q) \cup  (P\times \{e_2\}) \cup (\{e_1\} \times Q)\, .$

 \bigskip This means that if both $S$ and $T$ happen to be empty,  then $\lambda$ is  defined on the whole of $X$ and 
\eqref{eq:free} holds for every nontrivial  $k$ in $K$, as desired. 
 
 \bigskip {\it We assume from now on and until the end of Case 1 that $S$ is nonempty}. 
 
 \bigskip Consider $s\in S$. For each $j \in J$ we set  $\lambda(s,q_j^{-1})=1$. 
 
 \medskip Using \eqref{eq:piab} and \eqref{eq:ab}, we see that
\eqref{eq:free} will hold for

$$k= (s,q_j)\text{ and }k= (s,q_j^{-1})\text{ if }\ \lambda(s,q_j)U(s,q_j)\not\simeq U(s,q_j^{-1})\, ; $$
$$k= (s,e_2)\text{ if}\ \lambda(s,q_0)U(s,q_0)\not\simeq \big(\lambda(s,q_0)\, U(s,q_0)\big)\sp* \, .$$
For each $j \in J\, , \, j \neq 0$, we therefore pick
\begin{equation*}
\lambda(s,q_j)\in\T\setminus \Omega(s,q_j)\, .
\end{equation*}

\medskip
 We also pick $\lambda(s,q_0) \in\T\setminus(\Omega(s,q_0)\cup\Omega_1(s))$.  

\medskip Again, these choices are possible as all the involved $\Omega$'s are countable. 

\medskip Following this procedure for every $s \in S$, we achieve that
$\lambda$ is defined on $G'_1 \times Q$ in such a way that \eqref{eq:free} will hold for all $$k \in \big(G'_1\times ( \{e_2\} \cup Q)\big) \cup (\{e_1\} \times Q)\, .$$  
If $T$ happens to be empty, this means that $\lambda$ is defined on the whole of $X$ and 
  \eqref{eq:free} holds for every nontrivial  $k$ in $K$, as desired.

\bigskip {\it  Finally, we assume from now on and until the end of Case 1 that $T$ is also nonempty}.  

\medskip Consider $t \in T$. For each $i \in I$ we set  $\lambda(p_i^{-1}, t)=1$. 
 
 \medskip Using  \eqref{eq:piab} and \eqref{eq:ab}, we see that
\eqref{eq:free} will hold for
$$k=(p_i,t)\text{ and }k=(p_i^{-1},t)\text{ if }\ \lambda(p_i,t)U(p_i,t)\not\simeq U(p_i^{-1},t)\, ;$$
$$k=(e_1,t)\text{ if}\ \lambda(p_0,t)U(p_0,t)\not\simeq \big(\lambda(p_0,t) U(p_0,t)\big)\sp*\, .$$
For each $i \in I\, , \, i \neq 0$, we  pick
$\ \lambda(p_i,t)\in\T\setminus \Omega(p_i,t)\, .$

\medskip We also pick $\lambda(p_0,t) \in\T\setminus\big(\Omega(p_0,t)\cup\Omega_2(t)\big)$. 

\medskip Once again, these choices are possible as all the involved $\Omega$'s are countable. 

\medskip By doing this for every $t \in T$, we achieve that
$\lambda$ is defined on $(G'_1\times G'_2) \setminus (S\times T)$ and  \eqref{eq:free} will hold for all 

$$k \in \big(G'_1\times (\{e_2\} \cup Q)\big) \cup \big ((\{e_1\}\cup P) \times G'_2\big )\, .$$  
It remains  to define $\lambda$  on $S\times T$ in a way which ensures that \eqref{eq:free} also will hold for all
$k \in S \times T$.  

\medskip Let $t \in T$. We will below describe how to define $\lambda$  on $S\times \{t\}$ in a way which ensures that \eqref{eq:free} will hold for all $k \in S \times \{t\}$.  By following this procedure for each $t \in T$, the proof in Case 1 will then be finished.   

\medskip It is now appropriate to partition $S$ as $S = S' \sqcup S''$, where
$$S' = \{ s \in S \mid sp_0 \in P\}\, , \, S'' = \{ s \in S \mid sp_0 \in S\}\, . $$

\smallskip Assume that $s\in S'$. 

\medskip Using  \eqref{eq:piab} and \eqref{eq:abc}, we see that \eqref{eq:free} will hold for
\begin{equation*}
k= (s,t)\text{ if }\  \lambda(s,t) U(s,t) \big(\lambda(p_0,t)U(p_0,t)\big)\sp*\not\simeq  \lambda(sp_0,t)U(sp_0,t)\, .
\end{equation*}

Note here that $\lambda(sp_0,t)$ is already defined since $sp_0 \in P$. Further, as  $\lambda(sp_0,t)U(sp_0,t)$ is diagonalisable,  the set
\begin{equation*}
\Omega'(s,t)=\{\mu \in \T \mid \mu \, \big(\lambda(sp_0,t)U(sp_0,t)\big)\simeq U(s,t)\big(\lambda(p_0,t)U(p_0,t)\big)\sp* \, \}\, 
\end{equation*}
is countable. 
\medskip We can therefore pick  
\begin{equation*}
\lambda(s,t)\in \T\setminus \overline{\Omega'(s,t)}\, .
\end{equation*}
If $S'$ is nonempty, we can do this for each $s \in S'$ and $\lambda$ will then be defined on $S' \times \{t\}$ in such a way  that \eqref{eq:free} will hold for every $k \in S' \times \{t\}$.  

\bigskip
If $S''$ is empty, then $S'$ has to be nonempty and the proof of Case 1 is then finished. 

\medskip Assume now that 
$S''$ is nonempty and consider
 $s \in S''$, so $(sp_0)^2=e_1$. 
 
 \medskip One easily checks that this implies that $s\, {p_0}^n={p_0}^{-n}\, s$ for all $n \in \Z$. It is then almost immediate that 
 $S''(s) = \{s{p_0}^n \mid n \in \Z\}$ is a subset of $ S''$. 
 
 \medskip Furthermore, if $\tilde{s} \in S'' \setminus S''(s)$, then $S''(s)$ and $S''(\tilde{s})$ are disjoint. 
 
 \medskip Hence, as $S''$ is countable, we may pick a countable family $\{s_l\}_{l\in L}$ of distinct elements in $S''$ such that $S'' = \sqcup_{l\in L} S''(s_l)\, .$

\medskip Consider $l\in L$. To ease notation we write $s=s_l$. 

\smallskip We are going to define  $\lambda$ on $S''(s) \times \{ t\}$  in such a way that \eqref{eq:free} will hold for every $k \in  S''(s) \times \{ t\}$.
By doing this for each $l \in L$, $\lambda$ will then be defined on  $S'' \times \{t\}$ and \eqref{eq:free} will hold for every $k \in S'' \times \{t\}$.

Since $S\times \{t\} = (S'\times \{t\}) \sqcup (S'' \times \{t\})$, the proof of Case 1 will then be finished.

\bigskip

For each $n \in \Z$, using \eqref{eq:piab} and \eqref{eq:abc} (with $a=s{p_0}^n, b=t$ and $c=s{p_0}^{n\pm1}$), we see that  \eqref{eq:free} will hold for
\begin{equation*}
\begin{split}
k=(s{p_0}^n,t)\text{ if}\ &\lambda(s{p_0}^n,t) U(s{p_0}^n,t) \big(\lambda(s{p_0}^{n+1},t)U(s{p_0}^{n+1},t)\big)\sp* \not\simeq \lambda(p_0,t)U(p_0,t)\\
         \text{ or}\ & \lambda(s{p_0}^n,t) U(s{p_0}^n,t) \big(\lambda(s{p_0}^{n-1},t)U(s{p_0}^{n-1},t)\big)\sp* \not\simeq  \lambda(p_0^{-1},t)U(p_0^{-1},t) \, .
\end{split}
\end{equation*}

\medskip Suppose first that $p_0$ is aperiodic, so $S''(s) = \sqcup_{n \in \Z}\{ s{p_0}^n \}\,.$

\medskip 
We first set $\lambda(s,t)=1$. Then, for each $m\in \N$, we do inductively the following two steps:  

\medskip 
i) Define 
$$ \Omega^m(s,t)=\{\mu \in \T \mid \mu \, \big(\lambda(p_0,t)U(p_0,t)\big) \simeq\lambda(sp_0^{m-1},t)U(sp_0^{m-1},t) U(s{p_0}^m,t)\sp*\}$$
(which is countable) and pick $\lambda(s{p_0}^m,t)\in\T\setminus \Omega^{m}(s,t)$.

\medskip
ii) Define
$$ \Omega^{\m m}(s,t)=\{\mu \in \T \mid \mu \, \big(\lambda({p_0}^{-1},t)U({p_0}^{-1},t)\big)\simeq\lambda(s{p_0}^{\m m+1},t)U(s{p_0}^{\m m+1},t) \, U(s{p_0}^{\m m},t)\sp*\}$$
(which is countable) and pick $\lambda(s{p_0}^{\m m},t)\in\T\setminus \Omega^{\m m}(s,t)$.

\medskip 
Once this inductive process is finished, $\lambda$ is defined on $S''(s) \times \{t\}$ and we know that \eqref{eq:free} holds for every $k=(s{p_0}^{\pm(m-1)}, t), \, m\in \N$, i.e.\ for every $k \in  S''(s) \times \{ t \}$, as desired.

\bigskip Assume now that $p_0$ is periodic with period $N$. Note that $N\geq 3$ since $p_0 \in P$. 
The aperiodic case has to be modified as follows. 

\medskip Again, we first set $\lambda(s,t)=1$. Then, for each $m=1, \cdots, N-2$, we define inductively 
$$ \Omega^m(s,t)=\{\mu \in \T \mid \mu \, \big(\lambda(p_0,t)U(p_0,t)\big) \simeq\lambda(sp_0^{m-1},t)U(sp_0^{m-1},t) U(s{p_0}^m,t)\sp*\}$$
(which is countable)
and pick $\lambda(s{p_0}^m,t)\in\T \setminus \Omega^{m}(s,t)$. 

\medskip This ensures that \eqref{eq:free} holds for each $k=(s{p_0}^{m-1}, t), \, m=1, \cdots, N-2.$

\medskip We also define 
$$ \Omega^{N-1}(s,t)=\{\mu \in \T \mid \mu \, \big(\lambda(p_0,t)U(p_0,t)\big) \simeq\lambda(sp_0^{N-2},t)U(sp_0^{N-2},t) U(s{p_0}^{N-1},t)\sp*\}$$
(which is countable). If we pick  $\lambda(s{p_0}^{N-1},t)$  outside $ \Omega^{N-1}(s,t)$, then \eqref{eq:free}  will hold for $k=(s{p_0}^{N-2}, t)$. However, we want to pick  $\lambda(s{p_0}^{N-1},t)$ so that \eqref{eq:free} also holds  for $k=(s{p_0}^{N-1}, t)$.

\medskip Now, using \eqref{eq:piab} and \eqref{eq:abc} (with $a=s{p_0}^{N-1}, b=t$ and $c=s$), we see that  \eqref{eq:free} will hold for
$k=(s{p_0}^{N-1},t) \text{ if}\ $ 
$$\lambda(p_0,t) U(p_0,t) \not\simeq \lambda({p_0}^{N-1},t)U({sp_0}^{N-1},t) U(s,t)\sp* \, .$$
Hence we define 
\begin{equation*}
\Omega_N(s,t)=\{\mu \in \T \mid \mu \, \big(\lambda(p_0,t)U(p_0,t)\big) \simeq U(sp_0^{N-1},t)  U(s,t)\sp*\}
\end{equation*}
(which is countable) and pick 
\begin{equation*}
\lambda(sp_0^{N-1},t)\in\T\setminus \big(\Omega^{N-1}(s,t)\cup\overline{\Omega_N(s,t)}\, \big).
\end{equation*}
This choice does ensure that \eqref{eq:free} holds both for $k=(s{p_0}^{N-2}, t)$ and $k=(s{p_0}^{N-1}, t)$.	

\medskip
Hence,  $\lambda$ is defined on $S''(s) \times \{t\}$ and \eqref{eq:free} holds for every $k \in  S''(s) \times \{ t \}$. This finishes the proof of Case 1. 

\end{case}

\begin{case}
{\it  Either $P$ is nonempty and $Q$ is empty, or $P$ is empty and $Q$ is nonempty }.

\medskip Clearly, it suffices to consider the first alternative. We then pick $p_0 \in P$, $t_0\in T$ and set $x_0=[{p_0}^{-1}, t_0] \in X$. 
  
  \medskip We let   $\pi=\pi_{x_0}$ be a Choi representation of $A$ associated to $x_0$ and set $U(a,b) = U_{\pi}(a,b)$ for each $x=[a,b] \in X$. 

\medskip Our proof that $\lambda : X \to \T$ may be defined so that \eqref{eq:free} holds for each nontrivial $k \in K$ is quite similar to our proof of Case 1, but  some care is required and some repetitions seem unavoidable in our presentation.

\medskip 
For  $a \in G'_1\, ,\, t \in T\,$, we now set  
\begin{equation*}
\begin{split}
\Omega(a,t)&=\{\mu  \in \T \mid \mu \, U(a,t)\simeq U(a^{-1},t)\}\, ,\\
\Omega_1(a)&=\{\mu \in \T \mid \mu \, U(a,t_0)\simeq U(a^{-1},t_0)\sp*\, \} \, ,\\
\Omega_2(t)&=\{\mu \in \T \mid \mu \, U(p_0,t)\simeq \big(\mu \, U(p_0, t)\big)\sp*\, \} \, . \\
\end{split}
\end{equation*}

Note that if $(a, t)\neq (p_0^{-1}, t_0)$, then $\Omega(a,t)$ is countable. 
On the other hand,  $\Omega_1(a)$ is countable when $a\neq p_0^{-1}$, and 
 $\Omega_2(t)$ is always countable.

\medskip 
Let  $P= \sqcup_{i \in I} \, \{p_i, p_i^{\, -1}\}$ be an enumeration of $P$, where $I$ is 
a (finite or infinite) set of successive integers starting from 0.  

\medskip First, we set $\lambda({p_i}^{-1},t)=1$ for all $i\in I$ and $t \in T$.

\medskip Let $i \in I, t \in T$. Using  \eqref{eq:piab} and \eqref{eq:ab}, we see that \eqref{eq:free} will hold for
\begin{equation*}
\begin{split}
k= (p_i,t)\text{ and }k=({p_i}^{-1},t)\text{ if }\ &\lambda(p_i,t)U(p_i,t)\not\simeq U({p_i}^{-1},t)\, ;\\
k=(p_i,e_2)\text{ and }k=({p_i}^{-1},e_2)\text{ if }\ &\lambda(p_i,t_0)U(p_i,t_0)\not\simeq U({p_i}^{-1},t_0)\sp*\, ;\\
                         k=(e_1,t)\text{ if }\ &\lambda(p_0,t)U(p_0,t)\not\simeq\big(\lambda(p_0,t)U(p_0,t)\big)\sp*\, .\\            
\end{split}
\end{equation*}

Therefore, for each $i \in I, i\neq 0$, and $t \in T, t \neq t_0$, we pick
\begin{equation*}
\begin{split}
\lambda(p_i,t)&\in\T\setminus \Omega(p_i,t)\, ,\\
\lambda(p_i,t_0)&\in\T\setminus\big(\Omega(p_i,t_0)\cup\Omega_1(p_i)\big)\, ,\\
\lambda(p_0,t)&\in\T\setminus\big(\Omega(p_0,t)\cup\Omega_2(t)\big)\, .\\
\end{split}
\end{equation*}
Finally, we pick 

$\quad \quad\quad \quad\quad \quad \quad \quad \, \, \lambda(p_0,t_0)\in\T\setminus\big(\Omega(p_0,t_0)\cup\Omega_1(p_0)\cup\Omega_2(t_0)\big) \, .$

\medskip These choices ensure that $\lambda$ is defined on $P \times T \, $ and  \eqref{eq:free} will hold for 
all $k \in (P \times \big(T\cup \{e_2\})\big)  \cup (\{e_1\} \times T) $.

 \medskip This means that if $S$ happens to be empty,  $\lambda$ is  defined on the whole of $X$ and 
\eqref{eq:free} holds for every nontrivial  $k$ in $K$, as desired. 
 
 \bigskip {\it We assume from now on and until the end of Case 2 that $S$ is nonempty}.

 \medskip Consider $s\in S$. 
 Using \eqref{eq:piab} and \eqref{eq:ab}, we see that
\eqref{eq:free} will hold for

\medskip $ \quad \quad k=(s,e_2)\text{ if }\ \lambda(s,t_0)U(s,t_0)\not\simeq\big(\lambda(s,t_0) U(s,t_0)\big)\sp*\, .$

\medskip We will therefore pick $\lambda(s,t_0)$ in a subset of $\T\setminus \Omega_1(s)$. But which subset will depend on whether
$s$ belongs to $S'$ or $S''$, where $$S' = \{ s \in S \mid sp_0 \in P\} \, \, \text{and} \, \, S'' = \{ s \in S \mid sp_0 \in S\}\,  $$
(using the same notation as in Case 1).
 
\bigskip Assume that  $s\in S'$, $t \in T$.

\medskip As in Case 1, \eqref{eq:free} will hold for
\begin{equation*}
k= (s,t)\text{ if }\  \lambda(s,t) U(s,t) \big(\lambda(p_0,t)U(p_0,t)\big)\sp*\not\simeq  \lambda(sp_0,t)U(sp_0,t)\, .
\end{equation*}
 Again, we set \begin{equation*}
\Omega'(s,t)=\{\mu \in \T \mid \mu \, \big(\lambda(sp_0,t)U(sp_0,t)\big)\simeq U(s,t)(\lambda(p_0,t)U(p_0,t)\big)\sp* \, \}\, .
\end{equation*}
 
If $t=t_0$, then we pick  
$\lambda(s,t_0)\in \T\setminus \big( \Omega_1(s) \cup \overline{\Omega'(s,t_0)}\big)\, .$  

\smallskip Otherwise, we pick $\lambda(s,t)\in \T\setminus  \overline{\Omega'(s,t)}\, .$

\medskip If $S'$ is nonempty, we can do this for every $s \in S'$ and every $t \in T$. This ensures that $\lambda$ is defined on $S' \times T$ and that \eqref{eq:free} will hold for every $ k \in (S' \times (T \cup \{e_2\}) $. Hence, if $S''$ is empty, then $S'$ has to be nonempty and  the proof of Case 2 is finished.  

\bigskip Assume now that $S''$ is nonempty.  As in Case 1, we then pick a countable family $\{s_l\}_{l\in L}$ of distinct elements in $S''$ such that $S'' = \sqcup_{l\in L} \, S''(s_l)\, ,$ where  $S''(s) = \{s{p_0}^n \mid n \in \Z\}$ for $s \in S''$.

\medskip Consider $l\in L, \,  t \in T$ and set $s=s_l$. 

\medskip If $t=t_0$, then we pick  $\lambda(s,t_0)\in \T\setminus  \Omega_1(s)$. Otherwise, we set $\lambda(s,t)=1$.

\medskip Let $n \in \Z$. 
As in Case 1,   \eqref{eq:free} will hold for
\begin{equation*}
\begin{split}
k=(s{p_0}^n,t)\text{ if}\ &\lambda(s{p_0}^n,t) U(s{p_0}^n,t) \big(\lambda(s{p_0}^{n+1},t)U(s{p_0}^{n+1},t)\big)\sp* \not\simeq \lambda(p_0,t)U(p_0,t)\\
         \text{ or}\ & \lambda(s{p_0}^n,t) U(s{p_0}^n,t) \big(\lambda(s{p_0}^{n-1},t)U(s{p_0}^{n-1},t)\big)\sp* \not\simeq  \lambda(p_0^{-1},t)U(p_0^{-1},t) \, .
\end{split}
\end{equation*}

\medskip Suppose first that $p_0$ is aperiodic, so $S''(s) = \sqcup_{n \in \Z}\{ s{p_0}^n \}\,.$

\medskip 
Then, for each $m\in \N$, we proceed inductively and do the following two steps:

\medskip 
i) Define 
$$ \Omega^m(s,t)=\{\mu \in \T \mid \mu \, \big(\lambda(p_0,t)U(p_0,t)\big) \simeq\lambda(sp_0^{m-1},t)U(sp_0^{m-1},t) U(s{p_0}^m,t)\sp*\}\, .$$
If $t=t_0$, then we pick  $\lambda(s{p_0}^m,t_0)\in\T\setminus \big( \Omega_1(s{p_0}^m) \cup \Omega^{m}(s,t_0)\big)$. 

\smallskip Otherwise, we pick  $\lambda(s{p_0}^m,t)\in\T\setminus  \Omega^{m}(s,t)$.

\medskip
ii) Define
$$ \Omega^{\m m}(s,t)=\{\mu \in \T \mid \mu \, \big(\lambda({p_0}^{-1},t)U({p_0}^{-1},t)\big)\simeq\lambda(s{p_0}^{\m m+1},t)U(s{p_0}^{\m m+1},t) \, U(s{p_0}^{\m m},t)\sp*\}\, .$$
If $t=t_0$, then we pick  $\lambda(s{p_0}^{\m m},t_0)\in\T\setminus \big( \Omega_1(s{p_0}^{\m m}) \cup \Omega^{\m m}(s,t_0)\big)$.

\smallskip Otherwise, we pick $\lambda(s{p_0}^{\m m},t)\in\T\setminus \Omega^{\m m}(s,t)$.

\bigskip Assume next that $p_0$ is periodic with period $N$ $\geq 3$.

\medskip 
Then for each $m=1, \cdots, N-2$,  we proceed inductively and do the following:  we define
$$ \Omega^m(s,t)=\{\mu \in \T \mid \mu \, \big(\lambda(p_0,t)U(p_0,t)\big) \simeq\lambda(sp_0^{m-1},t)U(sp_0^{m-1},t) U(s{p_0}^m,t)\sp*\}\, .$$
If $t=t_0$, we pick $\lambda(s{p_0}^m,t_0)\in\T \setminus \big( \Omega_1(s{p_0}^m) \cup \Omega^{m}(s,t_0)\big)$. Otherwise,   
we pick  $\lambda(s{p_0}^m,t)\in\T \setminus \Omega^{m}(s,t)$.

\bigskip We also define 
$$ \Omega^{N-1}(s,t)=\{\mu \in \T \mid \mu \, \big(\lambda(p_0,t)U(p_0,t)\big) \simeq\lambda(sp_0^{N-2},t)U(sp_0^{N-2},t) U(s{p_0}^{N-1},t)\sp*\}\, .$$ 

\medskip 
As in Case 1,  \eqref{eq:free} will hold for
$k=(s{p_0}^{N-1},t) \text{ if}\ $ 
$$\lambda(p_0,t) U(p_0,t) \not\simeq \lambda({p_0}^{N-1},t)U({sp_0}^{N-1},t) U(s,t)\sp* \, .$$
So we define 
\begin{equation*}
\Omega_N(s,t)=\{\mu \in \T \mid \mu \, \big(\lambda(p_0,t)U(p_0,t)\big) \simeq U(sp_0^{N-1},t)  U(s,t)\sp*\} \, .
\end{equation*}
Now, if $t=t_0$, then we pick 
\begin{equation*}
\lambda(sp_0^{N-1},t_0)\in\T\setminus \big(\Omega(sp_0^{N-1}) \cup \Omega^{N-1}(s,t_0)\cup\overline{\Omega_N(s,t_0)}\big).
\end{equation*}
Otherwise,  we pick 
\begin{equation*}
\lambda(sp_0^{N-1},t)\in\T\setminus \big(\Omega^{N-1}(s,t)\cup\overline{\Omega_N(s,t)}\big).
\end{equation*}

\medskip 
Under both alternatives ($p_0$ being aperiodic or not),  these processes ensure that $\lambda$ is defined on $S''(s) \times \{t\}$ and that \eqref{eq:free} will hold for every $k \in  S''(s) \times (\{ t \} \cup \{ e_2\})$. 

\medskip After having done this for every $s=s_l \, (l \in L)$  and every $t \in T$, $\lambda$ is defined on $S'' \times T$ and  we know that \eqref{eq:free} will hold for every $k \in  S'' \times (T \cup \{ e_2\})$. 
 
\medskip Altogether, this means that $\lambda$ is defined on the whole of $G'_1 \times G'_2$ and  \eqref{eq:free}  holds for every nontrivial $k \in K$. This finishes the proof of Case 2.

\end{case}

\bigskip 
\begin{case}
{\it Both $P$ and $Q$ are empty}. 

\medskip This means that $G'_1= S$ and $G'_2=T$, i.e.\ all elements in both groups have order $2$, so $G_1$ and $G_2$ are abelian. Moreover, as one of them is assumed to have more than two elements, we may assume that $\lvert G_1\rvert \geq 4$ and $\lvert G_2\rvert  \geq 2$.

\medskip We pick $s_0 \in S$, $t_0\in T$ and set $x_0=[{s_0}, t_0] \in X$. 
  
  \medskip Next, we let   $\pi=\pi_{x_0}$ be a Choi representation of $A$ associated to $x_0$ and set $U(a,b) = U_{\pi}(a,b)$ for each $(a,b) \in S \times T = G'_1 \times G'_2$. 

\medskip
 Now, since $S$ is countable, it is not difficult to see that we may find a family $\{s_l\}_{l\in L}$  of distinct elements in $S \setminus \{s_0\} $ such that  
 $$S = \{s_0\}\, \sqcup \, \big( \, \sqcup_{l\in L} \{ s_l, s_0s_l\} \, \big)\, ,$$
where $L$ is a (finite or infinite) set of successive integers starting from 1.

\medskip
Let $t \in T$. 
Set $\lambda(s_0,t)=1$ and  $\lambda(s_l,t)=1$ for each $l\in L, \,  l\geq 2$.   

\medskip Using \eqref{eq:piab} and  \eqref{eq:ab}, we see that \eqref{eq:free} will  hold for
\begin{equation*}
\begin{split}
                    k= (e_1,t)\text{ if}\ &\lambda(s_1,t)U(s_1,t)\not\simeq(\lambda(s_1,t)U(s_1,t))\sp* \, .\\
\end{split}
\end{equation*}
Hence we set 
$\Omega(t)=\{\mu \in \T \mid \mu  \, U(s_1,t) \simeq (\mu \, U(s_1,t))\sp* \}$, which is countable, and  pick
\begin{equation*}
\lambda(s_1,t)\in\T\setminus\Omega(t).
\end{equation*}

\medskip  Consider now $l\in L$.  Using \eqref{eq:piab},  \eqref{eq:ab}, \eqref{eq:abc} and \eqref{eq:aecb}, we see that \eqref{eq:free} will  hold for
\begin{equation*}
\begin{split}
k=(s_0,t)\text{ and } k=(s_l,e_2)  \text{ if}\ &\lambda(s_0s_l,t)U(s_0s_l,t)\not\simeq U(s_0,t)(\lambda(s_l,t)U(s_l,t))\sp* \, ;\\
k=(s_0,e_2) \text{ and } k=(s_l,t) \text{ if}\ &  \lambda(s_0s_l,t)U(s_0s_l,t) \not\simeq  \lambda(s_l,t)U(s_l,t)U(s_0,t)\sp* \, ;\\
k=(s_0,t) \text{ and } k=(s_0s_l,e_2) \text{ if}\ &  \lambda(s_l,t)U(s_l,t)   \not\simeq U(s_0,t)(\lambda(s_0s_l,t)U(s_0s_l,t))\sp* \, ;\\
k=(s_0s_l,t)\text{ and } k=(s_0,e_2)\text{ if}\ &\lambda(s_l,t)U(s_l,t)\not\simeq\lambda(s_0s_l,t)U(s_0s_l,t)U(s_0,t)\sp* \, .
\end{split}
\end{equation*}

 For each $l \in L$, we therefore set 
\begin{equation*}
\begin{split}
\Omega_1(l,t)&=\{\mu \in \T \mid \mu \, U(s_ls_0,t) \simeq U(s_0,t)\big(\lambda(s_l,t)U(s_l,t)\big)\sp*\} \, ,\\
\Omega_2(l,t)&=\{\mu \in \T \mid \mu \, U(s_0s_l,t) \simeq  \lambda(s_l,t)U(s_l,t)U(s_0,t)\sp*\} \, ,\\
\Omega_3(l,t)&=\{\mu \in \T \mid   \mu \, \big(\lambda(s_l,t)U(s_l,t)\big) \simeq U(s_0,t)U(s_0s_l,t)\sp* \} \, ,\\
\Omega_4(l,t)&=\{\mu \in \T \mid \mu  \, \big(\lambda(s_l,t)U(s_l,t)\big) \simeq U(s_ls_0,t)U(s_0,t)\sp*\} \, .
\end{split}
\end{equation*}
All these sets are countable. Hence, for each $l \in L$, we can pick
\begin{equation*}
\lambda(s_0s_l,t)\in\T\setminus\Big(\Omega_1(l, t)\cup\Omega_2(l,t)\cup\Omega_3(l,t)\cup\overline{\Omega_4(l,t)}\,\Big)\, .
\end{equation*}
We have thereby defined $\lambda$ on $S \times \{t\}$ in such a way that  \eqref{eq:free} will hold for every $k \in \big(G_1 \times \{t\}\big) \, \sqcup \,  \big(S \times \{e_2\}\big)\, .$ By doing this for each $t \in T$,  $\lambda$ is defined on $S \times T = G'_1 \times G'_2$ and  \eqref{eq:free} holds for every nontrivial $k \in K$. This finishes the proof of Case 3 (and thereby the proofs of Proposition 1 and Theorem 2).

\end{case}

\end{proof}

\section{Some further aspects}

We believe that if $G$ is a countable group such that $C^*(G)$ is primitive, then $C^*(G)$ is antiliminary and has an uncountable family of  pairwise inequivalent,  irreducible faithful representations. It is not difficult to see that this true in the case where $G$ is nontrivial, icc and amenable (see below). As pointed out in \cite{BO}, this also holds when $G=\Z_2 * \Z_3$. The argument was based on the following observation, which goes back to the work of J. Glimm and J. Dixmier in the sixties. We recall that a representation of a C$^*$-algebra is called {\it essential} whenever its range contains no compact operators other than zero.

\bigskip 
{\bf Proposition 2}. \ 

\smallskip {\it Let $A$ be a primitive separable C$^*$-algebra and consider the set $ \ospec{A}=\{[\pi]\in\widehat{A}\mid\pi\text{ is faithful}\}$.  Then the following conditions are equivalent:
\begin{itemize}
\item[i)] $ | \,\ospec{A}\, | > 1$.
\item[ii)] Every faithful irreducible representation of $A$ is essential.
\item[iii)] $A$ has a faithful irreducible representation which is essential.
\item[iv)] $\ospec{A}$ is uncountable.
\end{itemize}

Moreover, if $A$ satisfies any of these conditions, then $A$ is antiliminary.}

\bigskip The implications $ii) \Rightarrow iii)$ and $iv) \Rightarrow i)$ are trivial. The implication $i) \Rightarrow ii)$ follows from \cite[Cor.\,\,4.1.10]{Di}, while $iii) \Rightarrow iv)$ follows from  \cite[Compl{\'e}ments\,\,4.7.2]{Di}. The final assertion follows from \cite[Compl\' ements 9.5.4]{Di}.

\bigskip For completeness we mention that  there is another way to show that a unital separable C$^*$-algebra is primitive and antiliminary. Indeed, using that primitivity and primeness are equivalent notions for separable C$^*$-algebras (see e.g. \cite{Ped}), one deduces that a separable unital C$^*$-algebra $A$ is primitive and antiliminary if and only if the  pure state space of $A$ is weak*-dense in the state space of  $A$ (cf.\ \cite[Lemme 11.2.4  and Compl\' ements 11.6.6]{Di}). H.\ Yoshizawa showed in \cite{Y} that the right side of this equivalence holds when $A=C^*(\F_2)$. 

\bigskip Now let $G=G_1*G_2$ be as in Theorem 2. It is conceivable that one might be able to check that condition i) in Proposition 2 holds for $A=C^*(G)$ by following the line of proof  used  in \cite{BO} when $G=\Z_2 * \Z_3$. However, in light of our proof of Theorem 2,  the necessary combinatorics  will certainly be very messy. 
We will instead use  the following well known lemma to check that condition ii) holds for $A=C^*(G)$ in many cases.

\bigskip
{\bf Lemma 1}.  {\it Let $A$ be a primitive, unital,  infinite-dimensional C$^*$-algebra. Assume that $A$ contains no nontrivial projections or that $A$ has a faithful tracial state. Then $A$ satisfies condition ii) in Proposition 2.} 

\bigskip
{\bf Proof}. 
For completeness, we give the proof. Let $\pi$ be a faithful irreducible representation of $A$ acting on a Hilbert space $\H$ and let $\K$ denote the compact operators on $\H$. Note that $\H$ is infinite-dimensional since $\pi(A)$ is infinite-dimensional.

\smallskip Assume first that $A$ contains no nontrivial projections. Since $\pi$ is faithful,  $\pi(A)$ contains no nontrivial projections. Hence $\pi(A) \cap \K = \{0\}$ (otherwise we would have $\K \subset \pi(A)$ by irreducibility, and $\pi(A$) would contain all finite-dimensional projections), so $\pi$ is essential.  

\smallskip Assume now that $A$ has a faithful tracial state $\tau$. Assume (for contradiction) that $\pi(A) \cap \K \neq \{0\}$. Then $\K \subset \pi(A)$. As is well known, when  $\H$ is infinite-dimensional, the only bounded trace on $\K$ is the zero map. Hence the restriction of $\tau$ to $\K$ must be zero. But $\K$ contains nontrivial projections and evaluation of $\tau$ on any of these  does not give zero since $\tau$ is faithful. This gives a contradiction, and it follows that $\pi$ is essential. 

\smallskip
\hfill$\square$

\bigskip
{\bf Corollary 1.} {\it Let $G=G_1*G_2$ satisfy the assumptions of Theorem 2. Assume also that $G_1$ and $G_2$ are both torsion-free.  Then $C^*(G)$ has no nontrivial projections. Moreover, it is antiliminary and has an uncountable family of of  pairwise inequivalent,  irreducible faithful representations}.

\bigskip
{\bf Proof}. The first assertion is mentioned by G.J.ÊMurphy  \cite[p. 703]{Mur}, where he refers to \cite{dHD}Ê and \cite{LP}Ê for a proof. It seems to us that this is somewhat unprecise\footnote{Here is an elaboration of this remark. R.\ Li and S.\ Pedersen introduce in \cite{LP} a certain property C for a countable group K which ensures that $C^*(K)$ has no nontrivial projections. Then they show that the free product of two countable groups with property C also has property C. However, we are not aware of any references showing that a torsion-free countable amenable group H has property C. In the paper by de la Harpe and Dykema \cite{dHD} that Murphy refers to, the proof that $ C^*(H) \simeq C_r^*(H)$  contains no nontrivial projections relies on the fact that countable amenable groups 
have the Haagerup property. As the Baum-Connes conjecture holds for any countable group with this property (as shown by Higson and Kasparov), it follows that $H$  also satisfies the Kadison-Kaplansky conjecture, i.e.\  $C_r^*(H)$ contains no nontrivial projections.}. We propose an alternative way to prove this assertion: 

\smallskip Since $G_1$ and $G_2$ are amenable,  $G$ has  the Haagerup property (\cite[Proposition 6.2.3]{CC}). Hence, as shown by N. Higson and G. Kasparov in \cite{HK},  $G$ satisfies the Baum-Connes conjecture. As $G$ is easily seen to be torsion-freee, $G$ also satisfies the Kadison-Kaplansky conjecture (see e.g.\ \cite{Va}), i.e.\  the reduced group C$^*$-algebra $C_r^*(G)$ contains no nontrivial projections. 

Moreover, as shown by J.L.\ Tu in \cite{Tu}, any group having the Haagerup property is K-amenable. It follows that the homomorphism $\lambda_{\ast}$ from $ K_0(C^*(G))$ to $ K_0(C_r^*(G))$  induced by the canonical map $\lambda : C^*(G) \to C_r^*(G)$ is an isomorphism. It is then straightforward to check that this implies that $C^*(G)$ has no nontrivial projections.   

\medskip Now, Theorem 2 says that $C^*(G)$ is primitive. The second assertion follows therefore from Proposition 2 in combination with the first assertion and Lemma 1. 

\smallskip
\hfill$\square$

\medskip To our knowledge, the class of countable discrete groups which are such that their full group C$^*$-algebras have a {\it faithful} tracial state has not been much studied. Clearly, it  does contain all countable amenable groups (as the full and the reduced group C$^*$-algebras agree for such groups, and the canonical tracial state on the reduced algebra is always faithful). Hence, if $H$ is nontrivial, icc and amenable, then $C^*(H)$ is primitive (cf.Ê\cite{Mur, Pa}) and Lemma 1  may be applied. Our assertion at the beginning of this section follows then from Proposition 2. On the other hand, this class also contains all free groups with countably many generators. This fact is due to Choi \cite[Corollary 9]{Cho} and may be put in a somewhat more general framework as follows.

\medskip We first recall that a C$^*$-algebra is called {\it residually finite-dimensional} (RFD) if it has a separating family of finite-dimensional representations (see e.g. \cite{EL}). Clearly, any abelian or finite-dimensional C$^*$-algebra is RFD. If $F$ is a free group  on countably many generators, then $C^*(F)$ is RFD (cf. \cite[Theorem 7]{Cho}). 
Moreover, the class of RFD C$^*$-algebras is closed under free products (see \cite[Theorem 3.2]{EL}). Finally, any unital RFD C$^*$-algebra has a faithful tracial state (see the proof of \cite[Corollary 9]{Cho}). Hence we get:

\bigskip 
{\bf Corollary 2.} {\it Consider $G=G_1*G_2$, where at least one of the $G_i$'s has more than two elements, and assume that $G_1$ (resp. $G_2$) is abelian or finite.  Then $C^*(G)$ is RFD, antiliminary and has an uncountable family of of  pairwise inequivalent,  irreducible faithful representations}.

\bigskip
{\bf Proof}.  It follows from Theorem 2 that $C^*(G)$ is primitive. Moreover, $C^*(G)=C^*(G_1)\ast C^*(G_2)$ is RFD since $C^*(G_1)$ and $C^*(G_2)$ are RFD. Hence $C^*(G)$ has a faithful tracial state, and the assertion follows from Proposition 2 combined with Lemma 1. 

\smallskip
\hfill$\square$

\end{flushleft}
\end{document}